\documentclass[12pt,a4paper]{amsart}
\usepackage{amsmath,amssymb,amsfonts}
\usepackage{mathpazo}
\usepackage[a4paper,margin=0.8in]{geometry}
\usepackage[breaklinks=true]{hyperref}
\usepackage{cleveref}
\allowdisplaybreaks[4]
\newtheorem{thm}{Theorem}[section]
\newtheorem*{thm*}{Theorem}
\newtheorem{prop}[thm]{Proposition}
\newtheorem{lem}[thm]{Lemma}
\newtheorem{cor}[thm]{Corollary}

\theoremstyle{definition}
\newtheorem{definition}[thm]{Definition}

\newtheorem{rem}[thm]{Remark}
\numberwithin{equation}{section}

\newcommand{\F}{\mathbb{F}}
\newcommand{\C}{\mathbb{C}}
\DeclareMathOperator{\tr}{tr}
\DeclareMathOperator{\Tr}{Tr}
\DeclareMathOperator{\rank}{rank}

\DeclareMathOperator{\GL}{GL}
\DeclareMathOperator{\Z}{\mathbb{Z}}
\DeclareMathOperator{\D}{\mathcal{D}}

\newcommand{\qbinom}[2]{\begin{bmatrix}#1\\#2\end{bmatrix}_q}

\begin{document}
\title{Spectra of gcd-graphs over a matrix ring}
\author{Tung T. Nguyen}

\address{Department of Mathematics, Elmhurst University, Illinois, USA}
 \email{tung.nguyen@elmhurst.edu}
  \keywords{Gcd graphs, Finite rings, Partial Traces, Character Sums}
\subjclass[2020]{Primary 05C50; Secondary 05C25, 15B33}
\date{}
\maketitle
\begin{abstract}
    In this short note, we calculate the spectra of gcd-graphs defined over a matrix ring. Using character sums, we show that this calculation can be reduced to certain counting problems. In particular, when the underlying ring is a matrix ring with coefficients in a finite field, we calculate these spectra explicitly. We will also show that two gcd-graphs over the same matrix ring are isospectral if and only if they have the same generating set. 
\end{abstract}

\section{Introduction} 
Gcd-graphs were first introduced by Klotz and Sander in \cite{klotz2007some}. These graphs lie at the intersection of algebra, number theory, and discrete analysis, and therefore they naturally exhibit many elegant properties. Ever since this article,  subsequent works have generalized the notion of a gcd-graph to an arbitrary finite ring $R$ (see \cite{minavc2024gcd, nguyengcd2026}). Let us first recall the theory of Cayley graphs over $R.$ Let $S$ be a symmetric subset of $(R,+)$ such that $0 \not \in S$. The Cayley graph $\Gamma(R,S)$ is defined as follows 

\begin{enumerate}
    \item The vertex set of $\Gamma(R,S)$ is $R.$
    \item Two vertices $u,v$ are adjacent if $u-v \in S.$
\end{enumerate}
We say that $\Gamma(R,S)$ is a gcd-graph if $S$ is stable under left and right actions of the unit group $R^{\times}$; i.e., $R^{\times} S R^{\times}=S.$
In this setting, the interplay between the additive and multiplicative structures of $R$ has played a fundamental role in understanding the geometry of their associated gcd-graphs.  Furthermore, when the underlying ring $R$ is a finite commutative Frobenius ring, we showed in previous work that the spectra of gcd-graphs over $R$ can be described by various arithmetical sums such as Gauss sums, Ramanujan sums, and Heilbronn sums (see \cite{nguyen2025supercharacters, nguyengcd2026}). The goal of this article is to extend these calculations to the general case where $R$ is a finite symmetric $\F_q$-algebra which could be non-commutative. When $R$ is a matrix ring with coefficients in a finite field, our approach gives complete and explicit formulas for the spectra of gcd-graphs defined over $R.$ To achieve this goal, we use the counting formula of Balasubramanian, Kaipa, and Khurana~\cite{BKK}, who also identify the relevant quantities with character sums over matrices of prescribed rank. Our approach to calculating the spectra is similar in spirit to the work of Eichler and Lamprecht; see the introduction of \cite{BKK} for a historical overview. Finally, we also show that the rank counts form a super collection and we use this fact to deduce that the degree of a gcd-graph over a fixed matrix ring determines its generating set. This proves a natural analogue of So's conjecture for gcd-graphs defined over a matrix ring.


\section{Spectra of a Cayley graph over a finite symmetric algebra}

In this section, we recall the theory of finite symmetric Frobenius algebras and the associated spectra of Cayley graphs defined over them. This theory is discussed in \cite{nguyen2026u}. 

Let $R$ be a finite-dimensional $\F_q$-algebra. We say that $R$ is a finite symmetric Frobenius $\F_q$-algebra if there exists a non-degenerate linear functional 
$\Psi:R\to\F_q$ such that  $\Psi(ab)=\Psi(ba).$ Here, non-degenerate means that the kernel of $\Psi$ does not contain any non-zero left ideal in $R.$ Equivalently, the map $R \to \text{Hom}_{\F_q}(R, \F_q)$ defined by $a \mapsto \Phi_a$ where $\Phi_a(u)=\Psi(au)$ is an isomorphism of $\F_q$-vector spaces. 
Let
\[
\widetilde\Psi=\Tr_{\F_q/\F_p}\circ\Psi:R\longrightarrow\F_p.
\]
Here $\Tr_{\F_q/\F_p}$ is the trace map from $\F_q$ to $\F_p.$ The pairing $(a,b)\mapsto\widetilde\Psi(ab)$ is also nondegenerate, showing that $R$ is also a finite symmetric $\F_p$-algebra. 

Let $\Gamma(R,S)$ be a Cayley graph over $R$, where $0\notin S$ and $S$ is stable under $\F_q^\times$. Let $\zeta_p=e^{2\pi i/p}$. The additive characters of $R$ are
\[
\chi_r(x)=\zeta_p^{\widetilde\Psi(rx)},\qquad r\in R.
\]
By the theory of $R$-circulant matrices (see \cite{circulant}), we know that the spectrum of $\Gamma(R,S)$ consists of the eigenvalues
\[
\lambda_r
=\sum_{s\in S}\zeta_p^{\widetilde\Psi(rs)}
=\sum_{s\in S}\zeta_p^{\Tr_{\F_q/\F_p}(\Psi(rs))},
\qquad r\in R,
\]
with one eigenvalue for each $r\in R$, counting multiplicities. We will show below that, in some cases, the calculation of $\lambda_r$ can be reduced to certain counting problems. To do so, we start with a simple lemma in character theory.

\begin{lem}\label{lem:character-sum}
Let $u\in\F_q$. Then
\[
\sum_{a\in\F_q^\times}\zeta_p^{\Tr_{\F_q/\F_p}(au)}
=\begin{cases}
q-1,&u=0,\\
-1,&u\ne0.
\end{cases}
\]
\end{lem}
\begin{proof}
The map
\[
\F_q\longrightarrow\C^\times,\qquad a\longmapsto\zeta_p^{\Tr_{\F_q/\F_p}(au)}
\]
is an additive character of $\F_q$. If $u=0$, it is the trivial character.
If $u\ne0$, it is nontrivial, so its sum over $\F_q$ is zero.
Removing the term $a=0$ gives the result.
\end{proof}

Let us introduce the counting numbers. 
\[
A_r=\#\{s\in S:\Psi(rs)=0\},\qquad
B_r=\#\{s\in S:\Psi(rs)\ne0\}.
\]
Then $A_r+B_r=|S|$. We are now ready to express $\lambda_r$ in terms of these counting numbers. 

\begin{prop}\label{prop:eigenvalue-count}
For each $r\in R$, we have
\[
\lambda_r=\frac{qA_r-|S|}{q-1}.
\]
\end{prop}
\begin{proof}
Let us write  $S/\F_q^\times=\{[s_1],\ldots,[s_d]\}$  for the scalar orbits. Since $0\notin S$, every orbit has $q-1$ elements.
Therefore
\begin{align*}
\lambda_r &=\sum_{i=1}^d\sum_{s\in[s_i]} \zeta_p^{\Tr_{\F_q/\F_p}(\Psi(rs))} =\sum_{i=1}^d\sum_{a\in\F_q^\times} \zeta_p^{\Tr_{\F_q/\F_p}(\Psi(ra s_i))}=\sum_{i=1}^d\sum_{a\in\F_q^\times} \zeta_p^{\Tr_{\F_q/\F_p}(a\Psi(rs_i))}.
\end{align*}
The condition $\Psi(rs)=0$ is constant on each scalar orbit. By Lemma~\ref{lem:character-sum},
\begin{align*}
\lambda_r
&=\frac{A_r}{q-1}(q-1)+\frac{B_r}{q-1}(-1) =A_r-\frac{B_r}{q-1} =A_r-\frac{|S|-A_r}{q-1}=\frac{qA_r-|S|}{q-1}.
\end{align*}
\end{proof}

\section{Applications to gcd-graphs over a matrix ring}
In this section, we focus on the case $R=M_n(\F_q).$ In this case, $M_n(\F_q)$ is a symmetric Frobenius ring with respect to the matrix trace map.  Let
\[
V_m=\{B\in M_n(\F_q):\rank(B)=m\} .\]

In \cite{NguyenNguyen}, using the theory of row and column operations, we show that $\Gamma(M_n(\F_q), S)$ is a gcd-graph if and only if 
 $S=\bigsqcup_{m\in\mathcal D}V_m,$ where $\mathcal D\subseteq\{1,\ldots,n\}$.  In this case, we will write $\Gamma(M_n(\F_q), \mathcal{D})$ for $\Gamma(M_n(\F_q), S).$

For $X\in M_n(\F_q)$, write
\[
A_X=\#\{B\in S:\tr(XB)=0\}.
\]

\begin{lem}\label{claim:rank}
If $\rank(X)=\rank(Y)$, then $A_X=A_Y$.
\end{lem}
\begin{proof}
We can write $X=PYQ$ for some $P,Q\in\GL_n(\F_q)$. Thus
\begin{align*}
A_X
&=\#\{B\in S:\tr(PYQB)=0\}\\
&=\#\{B\in S:\tr(YQBP)=0\}\\
&=\#\{C\in S:\tr(YC)=0\}\\
&=A_Y.
\end{align*}
Here $C=QBP$, and $B\mapsto QBP$ is a rank-preserving bijection of $S$.
\end{proof}

For $0\le k\le n$, let
\[
I_n(k)=\begin{pmatrix}I_k&0\\0&0\end{pmatrix}.
\]
By Lemma \ref{claim:rank}, we only need to calculate
\[
A_{S,k}:=A_{I_n(k)}
=\#\{B\in S:\tr(I_n(k)B)=0\}.
\]
If
\[
B=\begin{pmatrix}B_{11}&B_{12}\\B_{21}&B_{22}\end{pmatrix},
\qquad B_{11}\in M_k(\F_q),
\]
then
\[
I_n(k)B=\begin{pmatrix}B_{11}&B_{12}\\0&0\end{pmatrix},
\qquad
\tr(I_n(k)B)=\tr(B_{11})=\sum_{i=1}^k b_{ii}.
\]
Consequently,
\[
A_{S,k}=\sum_{m\in\mathcal D}A_{m,k},
\]
where
\begin{equation}\label{eq:Amk}
A_{m,k}
=\#\left\{B=(b_{ij})\in M_n(\F_q):
\rank(B)=m,\ \sum_{i=1}^k b_{ii}=0\right\}.
\end{equation}
Here $n$ and $q$ are fixed. In particular, $A_{m,0}=|V_m|$.

It turns out that there is a neat formula for $A_{m,k}$ recently discovered in~\cite{BKK}. Let us briefly recall this formula. First, let us write 
\[
a(N,r,q)=\prod_{j=0}^{r-1}\frac{(q^N-q^j)^2}{q^r-q^j},
\qquad 0\le r\le N,
\]
with $a(N,0,q)=1$. This is the number of rank-$r$ matrices in $M_N(\F_q)$. In particular, $|V_m|=a(n,m,q)$.

\begin{thm}[{\cite[Theorem 1.1]{BKK}}]\label{thm:partial-trace}
For $0\le m,k\le n$, we have
\[
A_{m,k}=\frac{a(n,m,q)+(q-1)g_{n,m,k}}{q},
\]
where
\[
g_{n,m,k}
=\sum_{i=\max\{0,m-(n-k)\}}^{\min\{m,k\}}
(-1)^i q^{\binom{i}{2}+k(m-i)}
\qbinom{k}{i}a(n-k,m-i,q).
\]
Here $\qbinom{k}{i}$ is the Gaussian binomial coefficient.
\end{thm}

The case $k=0$ follows directly from $g_{n,m,0}=a(n,m,q)$.

\begin{thm}\label{cor:spectrum}
The eigenvalues of the gcd-graph $\Gamma(M_n(\F_q),S)$ are
\[
\lambda_k=\sum_{m\in\mathcal D}g_{n,m,k},\qquad 0\le k\le n.
\]
For each $k$, the rank-$k$ characters contribute $a(n,k,q)$ copies of $\lambda_k$. If different values of $k$ give the same eigenvalue, their multiplicities are added.

\end{thm}
\begin{proof}
By Proposition~\ref{prop:eigenvalue-count} and Theorem~\ref{thm:partial-trace},
\begin{align*}
\lambda_k
&=\frac{qA_{S,k}-|S|}{q-1} = \sum_{m\in\mathcal D}
  \frac{qA_{m,k}-a(n,m,q)}{q-1} = \sum_{m\in\mathcal D}g_{n,m,k}.
\end{align*}
There is one character for each matrix $X\in M_n(\F_q)$, and there are $a(n,k,q)$ matrices of rank $k$.
\end{proof}
We have an immediate corollary. 
\begin{cor}
    A gcd-graph over $M_n(\F_q)$ is integral. 
\end{cor}

\begin{rem}
When $\mathcal D=\{n\}$, we recover the spectrum of the unitary Cayley graph over $M_n(\F_q)$. Indeed, setting $m=n$ in the formula for $g_{n,m,k}$, the only term in the sum is $i=k$. Thus
\begin{align*}
g_{n,n,k}
&=(-1)^k q^{\binom{k}{2}+k(n-k)}a(n-k,n-k,q)\\
&=(-1)^k q^{k(2n-k-1)/2}\gamma(n-k,q),
\end{align*}
where $\gamma(r,q)=|\GL_r(\F_q)|$ and $\gamma(0,q)=1$.  These are precisely the eigenvalues obtained in \cite[Theorem 1.1]{chen2022unitary}.  We remark that the proof presented in \cite{chen2022unitary} uses  an entrywise parametrization of the additive characters, which agrees with ours after transposing the indexing matrix. Their calculation therefore gives the same character sum $g_{n,n,k}$.
\end{rem}

\section{An analogue of So's conjecture}
Two graphs are called isospectral if their spectra are the same. Isomorphic graphs are isospectral but the converse is not true in general. 
For a gcd-graph over $\mathbb{Z}/N$, So conjectured in \cite{So2006} that two gcd-graphs $\Gamma(\Z/N, D_1)$ and $\Gamma(\Z/N, D_2)$ are isospectral if and only if $D_1=D_2.$ In this section, we will prove an analogue of So's conjecture for gcd-graphs over a matrix ring. More precisely, we will show the following. 

\begin{thm} \label{thm:isospectral-gcd-graphs}
    Suppose that $n \geq 1$. Let $\mathcal{D}_1, \mathcal{D}_2$ be two subsets of $\{1, 2, \ldots, n\}$ such that $\Gamma(M_n(\F_q), \D_1)$ and $\Gamma(M_n(\F_q), \mathcal{D}_2)$ are isospectral. Then $\D_1=\D_2.$
\end{thm}
In fact, we will show that having the same degree already guarantees that $\D_1 = \D_2.$  In order to prove this theorem, 
we first introduce the following definition which is motivated by our study of sedentariness on Cayley graphs \cite[Theorem 2.2]{MMN2026}; namely, if a walk-regular graph is super then it is sedentary at every vertex. 
\begin{rem}
After the completion of this paper, we learned that this concept has been studied in \cite[Definition~1.5]{MoniusSo2023}. We note, however, that our convention reverses the order of the super sequences defined in their paper. 
\end{rem}

\begin{definition}
A sequence of positive real numbers $a_1,\ldots,a_r$ is called \emph{super} if
\[
a_i>\sum_{j=i+1}^r a_j
\qquad\text{for each }1\leq i\leq r.
\]
A finite collection of positive real numbers is called super if its decreasing ordering is super. We say that a finite graph $G$ with spectrum $\{\lambda_1^{(m_1)}, \lambda_2^{(m_2)}, \ldots, \lambda_r^{(m_r)}\}$ is \emph{super} if the sequence $\{m_1, m_2, \ldots, m_r\}$ forms a super collection.
\end{definition}

For $0\leq k\leq n$, let $V_k$ be the set of matrices in $M_n(\F_q)$ that have rank $k$. Recall that $|V_0|=1$ and
\[
|V_k|=\prod_{j=0}^{k-1}\frac{(q^n-q^j)^2}{q^k-q^j}.
\]
We will show that these numbers form a super collection.

\begin{prop}\label{thm:super-rank-counts}
Suppose that $n\geq2$. Then the collection $\{|V_0|,\ldots,|V_n|\}$ is super. More precisely:
\begin{enumerate}
\item If $q\geq3$, then
\[
|V_n|,|V_{n-1}|,\ldots,|V_0|
\]
is a super sequence.
\item If $q=2$, then
\[
|V_{n-1}|,|V_n|,|V_{n-2}|,\ldots,|V_0|
\]
is a super sequence.
\end{enumerate}
\end{prop}

\begin{proof}
Write $N_k=|V_k|$. By the rank-counting formula, we have
\[
\frac{N_k}{N_{k-1}}=\frac{q^{k-1}(q^{n-k+1}-1)^2}{q^k-1}.
\]
Consequently, for $1\leq k\leq n-1$,
\[
\frac{N_k}{N_{k-1}}>\frac{(q^{n-k+1}-1)^2}{q}\geq\frac{(q^2-1)^2}{q}\geq\frac92>2.
\]
It follows that
\begin{equation}\label{eq:lower-rank-dominance}
\sum_{j=0}^{k-1}N_j<N_k\sum_{\ell=1}^k2^{-\ell}<N_k
\qquad(1\leq k\leq n-1).
\end{equation}

Let us focus on the first case, where $q\geq3$. Using the  inequality
\[
\prod_{i=1}^n(1-x_i)\geq1-\sum_{i=1}^n x_i
\qquad(0\leq x_i\leq1),
\]
we obtain
\begin{align*}
\frac{N_n}{|M_n(\F_q)|}
&=\prod_{i=1}^n(1-q^{-i}) \geq1-\sum_{i=1}^n q^{-i} \geq 1-\frac{1-q^{-n}}{q-1}\\
&>1-\frac1{q-1}\geq\frac12.
\end{align*}
Since $\sum_{k=0}^nN_k=|M_n(\F_q)|$, this gives $N_n>\sum_{k=0}^{n-1}N_k.$  Together with \eqref{eq:lower-rank-dominance}, this proves the first assertion.

Let us prove the second assertion, where $q=2$. For each integer $r\geq2$, define
\[
b_r=\frac{\#\{A\in M_r(\F_2):\rank(A)=r-1\}}{2^{r^2}}.
\]
By the rank-counting formula,
\[
b_r=(1-2^{-r})^2\prod_{i=2}^{r-1}(1-2^{-i}).
\]
Therefore,
\[
\frac{b_{r+1}}{b_r}
=\frac{(1-2^{-(r+1)})^2}{1-2^{-r}}
=1+\frac{2^{-2r-2}}{1-2^{-r}}>1.
\]
Thus the sequence $(b_r)_{r\geq2}$ is strictly increasing. Applying this to our fixed matrix size $n$, we obtain
\[
\frac{N_{n-1}}{2^{n^2}}
=b_n\geq b_2
=\left(1-\frac14\right)^2
=\frac9{16}>\frac12.
\]
Since $\sum_{j=0}^nN_j=2^{n^2}$, it follows that
\[
N_{n-1}>N_n+\sum_{j=0}^{n-2}N_j.
\]

It remains to compare $N_n$ with the lower-rank counts. The ratio formula gives
\[
\frac{N_n}{N_{n-2}}=\frac{9\cdot2^{2n-3}}{(2^n-1)(2^{n-1}-1)}>\frac94>2.
\]
On the other hand, the inequalities $N_k>2N_{k-1}$ imply
\[
\sum_{j=0}^{n-2}N_j\leq N_{n-2}\sum_{\ell=0}^{n-2}2^{-\ell}<2N_{n-2}<N_n.
\]
Together with \eqref{eq:lower-rank-dominance}, these inequalities prove that $N_{n-1},N_n,N_{n-2},\ldots,N_0$ is super.
\end{proof}

\begin{rem}
When $n=1$, we have $|V_0|=1$ and $|V_1|=q-1$. Thus the rank counts are also super when $q>2$. The only exception is $(n,q)=(1,2)$, where both counts equal $1$.
\end{rem}

\begin{cor}
Every gcd-graph over $M_n(\F_q)$ is super, provided $(n,q)\ne(1,2)$.
\end{cor}

\begin{proof}
For each $k$, the characters indexed by rank-$k$ matrices contribute $|V_k|$ copies of the corresponding eigenvalue. If several ranks give the same eigenvalue, their contributions are added. Thus the multiplicities of the distinct eigenvalues are obtained by merging terms of the super rank-count sequence.

We claim that merging terms preserves the super property. Indeed, let $a_1,\ldots,a_r$ be super, and partition $\{1,\ldots,r\}$ into nonempty blocks $B_1,\ldots,B_t$, ordered so that
\[
\min B_1<\cdots<\min B_t.
\]
Set $b_j=\sum_{i\in B_j}a_i$. If $s=\min B_j$, then
\[
b_j\geq a_s>\sum_{i>s}a_i\geq\sum_{\ell>j}b_\ell.
\]
Therefore $b_1,\ldots,b_t$ is super. Applying this observation to the rank counts proves the result.
\end{proof}
 We are now ready to prove Theorem \ref{thm:isospectral-gcd-graphs}.
\begin{proof}[Proof of Theorem~\ref{thm:isospectral-gcd-graphs}]
 Suppose that $\Gamma(M_n(\F_q), \D_1)$ and $\Gamma(M_n(\F_q), \D_2)$ are isospectral. Since each graph is regular, its largest adjacency eigenvalue is its degree. Therefore,
\[
\sum_{m\in D_1}|V_m|=\sum_{m\in D_2}|V_m|.
\]
Suppose first that $(n,q)\ne(1,2)$. We have shown that the rank counts $|V_0|,\ldots,|V_n|$ form a super collection. We claim that this implies $\D_1=\D_2$. Indeed, suppose to the contrary that $\D_1 \neq \D_2$.  Let us choose $k\in D_1\triangle D_2$ such that $|V_k|$ is largest among the rank counts indexed by $\D_1\triangle \D_2$. Interchanging $\D_1$ and $\D_2$ if necessary, we may assume that $k\in \D_1\setminus \D_2$. Cancelling the common terms in the degree equality gives
\[
\sum_{m\in \D_1\setminus \D_2}|V_m|
=
\sum_{m\in \D_2\setminus \D_1}|V_m|.
\]
However, by the super property and the choice of $k$,
\[
\sum_{m\in \D_1\setminus \D_2}|V_m|
\geq |V_k|
>
\sum_{m\in \D_2\setminus \D_1}|V_m|,
\]
a contradiction. Hence $\D_1=\D_2$.

Finally, when $(n,q)=(1,2)$, the only possible rank lists are $\varnothing$ and $\{1\}$, which give graphs of degrees $0$ and $1$, respectively. Thus the conclusion holds in this case as well.
\end{proof}

\section*{Acknowledgements}
The author used ChatGPT to assist with locating references, reviewing the exposition, and running computational checks. He thanks Dr. Nguyen Duy T{\^a}n for enlightening discussions around this topic. This project is supported by a Small Research Grant at Elmhurst University.

\end{document}